\documentstyle[amscd,amssymb,verbatim,epsf,graphics]{amsart}
\newcommand{\R}{\mathbb{R}}
\newcommand{\N}{\mathbb{N}}
\begin{document}

\theoremstyle{plain}
\newtheorem{Thm}{Theorem}
\newtheorem{Cor}{Corollary}
\newtheorem{Con}{Conjecture}
\newtheorem{Main}{Main Theorem}
\newtheorem{Exmp}{Example}
\newtheorem{Lem}{Lemma}
\newtheorem{Prop}{Proposition}

\theoremstyle{definition}
\newtheorem{Def}{Definition}
\newtheorem{Note}{Note}

\theoremstyle{remark}
\newtheorem{notation}{Notation}
\renewcommand{\thenotation}{}

\errorcontextlines=0
\numberwithin{equation}{section}
\renewcommand{\rm}{\normalshape}%

\title[A study on downward half Cauchy sequences]%
{   A study on downward half Cauchy sequences}
\author[Huseyin Cakalli, Maltepe University, Istanbul-Turkey]{Huseyin Cakalli\\
          Maltepe University, Maltepe, Istanbul-Turkey}
\address{Huseyin Cakalli\
          Maltepe University, Marmara E\u{g}\.{I}t\.{I}m K\"oy\"u, TR 34857, Maltepe, \.{I}stanbul-Turkey Phone:(+90216)6261050 ext:2311, fax:(+90216)6261113
}
\email{huseyincakalli@@maltepe.edu.tr; hcakalli@@gmail.com}

\date{\today}

\begin{abstract}
In this paper, we introduce and investigate the concepts of down continuity and down compactness. A real valued function $f$ on a subset $E$ of $\R$, the set of real numbers is down continuous if it preserves downward half Cauchy sequences, i.e.  the sequence $(f(\alpha_{n}))$ is downward half Cauchy whenever $(\alpha_{n})$ is a downward half Cauchy sequence of points in $E$, where a sequence $(\alpha_{ k})$ of points in $\R$ is called downward half Cauchy if for every $\varepsilon>0$ there exists an $n_{0}\in{\N}$ such that $\alpha_{m}-\alpha_{n} <\varepsilon$ for  $m \geq n \geq n_0$. It turns out that the set of down continuous functions is a proper subset of the set of continuous functions.
\end{abstract}

\maketitle

\section{Introduction}
\normalfont
The concept of continuity and any concept involving continuity play a very important role not only in pure mathematics, but also in other branches of sciences involving mathematics especially in computer science, dynamical systems, economics, information theory, biological science.

Using the idea of continuity of a real function in terms of sequences, many kinds of continuities were introduced and investigated, not all but some of them we recall in the following: slowly oscillating continuity  (\cite{CakalliSlowlyoscillatingcontinuity}, \cite{VallinCreatingslowlyoscillatingsequencesandslowlyoscillatingcontinuousfunctions}), quasi-slowly oscillating continuity (\cite{CanakandDik}), ward continuity (\cite{CakalliForwardcompactness}, \cite{CakalliForwardcontinuity},\cite{BurtonandColemanQuasiCauchysequences}), $p$-ward continuity (\cite{CakalliVariationsonquasiCauchysequences}), $\delta$-ward continuity (\cite{CakalliDeltaquasiCauchysequences}),  $\delta^{2}$-ward continuity (\cite{BrahaandCakalliAnewtypecontinuityforrealfunctions},  statistical ward continuity (\cite{CakalliStatisticalwardcontinuity}, \cite{CakalliStatisticalquasiCauchysequences}), $\lambda$-statistical ward continuity ((\cite{CakalliandSonmezandArasLambdastatisticallywardcontinuity}, \cite{MursaleenLamdastatisticalconvergence}),
$\rho$-statistical ward continuity (\cite{CakalliAVariationonStatisticalWardContinuity}), arithmetic continuity (\cite{CakalliAvariationonarithmeticcontinuity}) strongly lacunary ward continuity (\cite{CakalliNthetawardcontinuity}, \cite{CakalliandKaplanAstudyonNthetaquasiCauchysequences}, \cite{KaplanandCakalliVariationsonstronglylacunaryquasiCauchysequencesAIP}, \cite{KaplanandCakalliVariationsonstronglacunaryquasiCauchysequencesJNonlinearSciAppl}, and
\cite{KaplanandCakalliVariationsonstronglylacunaryquasiCauchysequencesAIP}), lacunary statistical ward continuity (\cite{CakalliandArasandSonmezLacunarystatisticalwardcontinuity},  \cite{CakalliandKaplanAvariationonlacunarystatisticalquasiCauchysequencesCommunications}, and  \cite{YildizIstatistikselboslukludelta2quasiCauchydizileri}), downward statistical continuity (\cite{CakalliUpwardanddownwardstatisticalcontinuities}), lacunary statistical downward continuity (\cite{CakalliandMucukLacunarystatisticallyupwardanddownwardhalfquasiCauchysequencesJofAnalysis}) which enabled some authors to obtain conditions on the domain of a function to be uniformly  continuous in terms of sequences in the sense that a function preserves a certain kind of sequences (see for example \cite[Theorem 6]{VallinCreatingslowlyoscillatingsequencesandslowlyoscillatingcontinuousfunctions},\cite[Theorem 1 and Theorem 2]{BurtonandColemanQuasiCauchysequences},\cite[Theorem 2.3]{CanakandDik}.

The purpose of this paper is to introduce the concepts of down continuity of a real function, and the concept of down compactness of a subset of $R$   which cannot be given by means of a sequential method $G$ and prove interesting theorems.

\section{Down compact sets}
Throughout this paper, $\mathbb{N}$ and $\R$ will denote the set of positive integers, and the set of real numbers, respectively.
Modifying the definition of a \textit{downward Cauchy sequence} introduced in \cite{ReillyandSubrahmanyamandVamanamurthyCauchysequencesinquasipseudometricspaces} and  \cite{CollinsandZimmerAnasymmetricArzelaAscolitheorem}, we state a  definition of downward half Cauchyness  of a real sequence given in \cite{PalladinoOnhalfcauchysequences} in the following.

\begin{Def}A sequence $(\alpha_{k})$ of points in $\R$ is called \textit{downward half Cauchy} if for each $\varepsilon>0$ there exists an $n_{0}\in{\mathbb{N}}$ so that $\alpha_{m}-\alpha_{n} <\varepsilon$ for $m \geq n \geq n_0$.
\end{Def}
We note that a sequence $(\alpha_{k})$ is downward half Cauchy if and only if $(-\alpha_{k})$ is upward half Cauchy. Trivially, any Cauchy sequence is downward  half Cauchy, so is a convergent sequence, but there are downward half Cauchy sequences which are not Cauchy. For example, the sequence $(a_{k} )=(-k)$ is  downward half Cauchy, but not Cauchy. Thus the set of downward half Cauchy sequences is a proper subset of the set of Cauchy sequences. Any subsequence of a  downward half Cauchy sequence is downward half Cauchy. We note that the sum of two downward  half Cauchy sequences is downward half Cauchy, the product of two bounded downward half Cauchy sequences is downward half Cauchy. For any positive real number $c$ and for any downward half Cauchy sequence $(\alpha_{k})$, the sequence $(c\alpha_{k})$ is downward  half Cauchy.

\begin{Def}
A subset $A$ of $\R$ is called down compact if any sequence of points in $E$ has a downward half Cauchy subsequence.
\end{Def}
First, we note that any finite subset of $\R$ is down   compact, the union of two down compact subsets of $\R$ is down compact, the intersection of any family of down compact subsets of $\R$ is down   compact,  and any subset of a down   compact set is down compact. Any statistically ward compact subset of $\R$ is down  compact. An  $I$-sequentially compact subset of $\R$ is down compact for a nontrivial admissible ideal $I$ (see \cite{CakalliandHazarikaIdealquasiCauchysequences}, \cite{SavasandCakalliIdealstatisticalquasiCauchysequences}, and \cite{CakalliAvariationonwardcontinuity}). Furthermore any bounded subset of $\R$ is down compact, any slowly oscillating compact subset of $\R$ is down compact (see \cite{CakalliSlowlyoscillatingcontinuity} for the definition of slowly oscillating compactness), any above bounded subset of $\R$ is down compact. These observations suggest to us the following.
\begin{Thm} \label{Theodowncompactiffboundedabove}
A subset of $\R$ is down compact if and only if it is bounded above.
\end{Thm}
\begin{pf}
The proof can be obtained by using a contradiction method.

\end{pf}

\section{Down continuous functions}
A real valued function $f$ defined on a subset of $\R$ is continuous if and only if, for each point $\ell$ in the domain, $\lim_{n\rightarrow\infty}f(\alpha_{n})=f(\ell)$ whenever $\lim_{n\rightarrow\infty}\alpha_{n}=\ell$. This is equivalent to the statement that $(f(\alpha_{n}))$ is a convergent sequence whenever $(\alpha_{n})$ is. This is also equivalent to the statement that $(f(\alpha_{n}))$ is a Cauchy sequence whenever $(\alpha_{n})$ is Cauchy provided that the domain of the function is closed. These known results for continuity for real functions in terms of sequences might suggest to us introducing a new type of continuity, namely, down continuity.
\begin{Def}
A function $f: E\rightarrow \mathbb{R}$ is called \textit{down continuous} on a subset of $\R$, if it preserves downward half Cauchy sequences, i.e.  the sequence $(f(\alpha_{k}))$ is downward half Cauchy whenever $(\alpha_{k})$ is a downward half Cauchy sequence of points in $E$.
\end{Def}
It should be noted that down continuity cannot be given by any $G$-continuity in the manner of \cite{ConnorandGrosseErdmannSequentialdefinitionsofcontinuityforrealfunctions}. We see that the composition of two down continuous functions is down continuous, and for every positive real number $c$, $cf$ is down continuous, whenever $f$ is down  continuous.\\
We see in the following that the sum of two down continuous functions is down continuous.
\begin{Thm} \label{TheoThesumofdowncontinuousfunctionsisdowncontinuous}
If $f$ and $g$ are down continuous functions, then $f+g$ is down continuous.
\end{Thm}
\begin{pf}
The proof can be obtained easily, so is omitted.
\end{pf}

The case for the product of functions is different. If $f$ and $g$ are bounded positive real valued down continuous functions, then the product of $f$ and $g$ is down continuous.

In connection with downward half Cauchy sequences, and convergent sequences the problem arises to investigate the following types of  "continuity" of functions on $\R$:
\begin{description}
\item[($\delta ^{-} $)] $(\alpha_{n}) \in {\Delta ^{-} } \Rightarrow (f(\alpha_{n})) \in {\Delta ^{-} }$
\item[($\delta ^{-}  c$)] $(\alpha_{n}) \in {\Delta ^{-}} \Rightarrow (f(\alpha_{n})) \in {c}$
\item[$(c\delta ^{-} )$] $(\alpha_{n}) \in {c} \Rightarrow (f(\alpha_{n})) \in {\Delta ^{-}}$.
\end{description}
We see that $(c)$ can be replaced by not only $\rho$-statistical continuity, but also lacunary statistical continuity, strongly lacunary-sequential continuity, $I$-sequential continuity, and more generally  $G$-sequential continuity (see \cite{CakalliSequentialdefinitionsofcompactness, CakalliOnGcontinuity}).
We see that $(\delta ^{-} )$ is down continuity of $f$. It is easy to see that $(\delta ^{-}  c)$ implies $(\delta ^{-} )$; $(\delta ^{-} )$ does not imply $(\delta ^{-}  c)$; $(\delta ^{-} )$ implies $(c\delta ^{-} )$; $(c\delta ^{-} )$ does not imply $(\delta ^{-} )$; $(\delta ^{-}  c)$ implies $(c)$, and $(c)$ does not imply $(\delta ^{-}  c)$; and $(c)$ implies $(c\delta ^{-})$. \\
Now we give the implication $(\delta ^{-} )$ implies $(c)$, i.e. any down continuous function is continuous.
\begin{Thm} \label{TheoDowncontinuityimpliesordinarycontinuity} Any down continuous function is continuous.
\end{Thm}
\begin{pf}
The proof can be obtained easily, so is omitted.
\end{pf}
We see in the following example that the converse of the preceding theorem is not always true.
\begin{Exmp}
The continuous function $f$ defined by $f(x)=-x$ for every $x \in {\mathbb{R}}$ is not down continuous.
\end{Exmp}

\begin{Thm} \label{TheoDowncontinousimageofdownwardcompactsubset} Down continuous image of any down half compact subset of $\R$ is down half compact.
\end{Thm}
\begin{pf}
Let $E$ be a subset of $\R$, $f:E\longrightarrow$ $\R$ be a down continuous function, and $A$ be a down half compact subset of $E$. Take any sequence $\boldsymbol{\beta}=(\beta_{n})$ of points in $f(A)$. Write $\beta_{n}=f(\alpha_{n})$, where $\alpha_{n}\in {A}$ for each $n \in{\N}$, $\boldsymbol{\alpha}=(\alpha_{n})$. Down half compactness of $A$ implies that there is a downward half Cauchy subsequence $\boldsymbol{\xi}$ of the sequence of $\boldsymbol{\alpha}$. Write $\boldsymbol{\eta}=(\eta_{k})=f(\boldsymbol{\xi})=(f(\xi_{k}))$. Then $\boldsymbol{\eta}$ is a down half Cauchy subsequence of the sequence $f(\boldsymbol{\beta})$. This completes the proof of the theorem.
\end{pf}
We note that down continuous image of any $N_{\theta}$-sequentially compact subset of $\R$ is $N_{\theta}$-sequentially compact, and
down continuous image of any $\rho$-statistically sequentially compact subset of $\R$ is $\rho$-statistically sequentially compact. Furthermore down continuous image of any $G$-statistically sequentially connected subset of $\R$ is $G$-sequentially connected (see \cite{CakalliSequentialdefinitionsofconnectedness}, \cite{CakalliandMucukConnectednessviaasequentialmethod}, \cite{MucukandCakalliGsequentiallyconnectednessfortopologicalgroupswithoperations}, and \cite{MucukandSahanOnGsequentialcontinuity}).
On the other hand, uniform continuity does not imply down continuity. The function $f: \mathbb{R} \rightarrow \mathbb{R}$ defined by $f(x)=-x$ for every $x\in{\mathbb{R}}$ is uniformly continuous, but not down continuous.

Now we have the following result related to uniform convergence, namely, the uniform limit of a sequence of down continuous functions is down continuous.
\begin{Thm} \label{TheoUniformlimitofdowncont} If $(f_{n})$ is a sequence of down continuous functions defined on a subset $E$ of $\R$ and $(f_{n})$ is uniformly convergent to a function $f$, then $f$ is down continuous on $E$.
\end{Thm}
\begin{pf}
The proof can be obtained directly so is omitted.
\end{pf}

\maketitle
\section{Conclusion}
In this paper, mainly a new types of continuity, namely down continuity of a real function, and down compactness of a subset of $\R$ are introduced and investigated. In this investigation we have obtained results related to down continuity, some other kinds of continuities via downward half Cauchy sequences, convergent sequences, statistical convergent sequences, lacunary statistical convergent sequences of points in $\R$. We note that the set of down continuous functions is a proper subset of the set of ordinary continuous functions. The term downward half Cauchy sequence can be considered to be associated with above boundedness of the underlying space, whereas the term Cauchy sequence is traditionally associated with the completeness of the underlying space. Instead of downward half Cauchy sequence, one might tend to consider an asymptotically non-decreasing sequence in the sense that a sequence $(\alpha_{n})$ is asymptotically non-decreasing if $\liminf_{k} \inf \{\alpha_{m} -\alpha_{n} : m > n\geq k, k\in{\N}\}\geq 0$, and investigate functions that preserve asymptotically non-decreasing sequences. As the set of downward half Cauchy sequences is different from the set of asymptotically non-decreasing sequences, a further new research would be carried out if one allows to find new results. We suggest to investigate downward half Cauchy sequences of fuzzy points in asymmetric fuzzy spaces (see \cite{CakalliandDasFuzzycompactnessviasummability}, \cite{ArasandSonmezandCakalliAnapproachtosoftfunctions}, \cite{KocinacSelectionpropertiesinfuzzymetricspaces} , \cite{ErdemandArasandCakalliandSonmezSoftmatricesonsoftmultisetsinanoptimaldecisionprocess} for the definitions and related concepts in fuzzy setting and soft setting). We also suggest to investigate downward half Cauchy double sequences (see for example \cite{CakalliandPattersonFunctionspreservingslowlyoscillatingdoublesequences}, \cite{PattersonandCakalliQuasiCauchydoublesequences}, \cite{PattersonandSavasRateofPconvergenceoverequivalenceclassesofdoublesequencespaces},   \cite{DjurcicandKocinacandZizovicDoublesequencesandselections} and \cite{PattersonandSavasAsymptoticEquivalenceofDoubleSequences} for the definitions and related concepts in the double sequences case). For another further study, we suggest to investigate downward half Cauchy sequences of points in an asymmetric cone metric space since in a cone metric space the notion of a downward half Cauchy sequence coincides with the notion of a Cauchy sequence, and therefore  down continuity coincides with ordinary continuity on complete subsets of $\R$ (see \cite{PalandSavasandCakalliIconvergenceonconemetricspaces}, \cite{CakalliandSonmezandGenc}, \cite{SonmezParacompactness},  \cite{SonmezandCakalliConenormedspacesandweightedmeans}, \cite{CakalliOnDeltaquasislowlyoscillatingsequences}, \cite{YayingandHazarikandCakalliNewresultsinquasiconemetricspaces}, and \cite{CakalliandSonmezSlowlyoscillatingcontinuityinabstractmetricspaces}).


\begin{thebibliography}{99}
\bibitem{ArasandSonmezandCakalliAnapproachtosoftfunctions} C.G. Aras, A. Sonmez, H. Cakalli, {\em An approach to soft functions}, J. Math. Anal. \textbf{8}, 2, 129-138, (2017).
\bibitem{BrahaandCakalliAnewtypecontinuityforrealfunctions} Naim L. Braha, H. Cakalli, {\em A new type continuity for real functions}, J. Math. Anal. \textbf{7},  6, 68-76, (2016).
\bibitem{CollinsandZimmerAnasymmetricArzelaAscolitheorem} J. Collins and J. Zimmer, An asymmetric Arzelà-Ascoli theorem, Topology and its Applications \textbf{154} (2007), 2312-2322.
\bibitem{ConnorandGrosseErdmannSequentialdefinitionsofcontinuityforrealfunctions} J. Connor  and K.G. Grosse-Erdmann, Sequential definitions of continuity for real functions, Rocky Mountain J. Math., \textbf{33}(1) (2003), 93-121.
\bibitem{CakalliForwardcompactness} \c{C}akalli H. Forward compactness, Conference on Summability and Applications, Shawnee State University, November 6 November 8, 2009. http://webpages.math.luc.edu/~mgb/ShawneeConference/Articles/HuseyinCakalliOhio.pdf.
\bibitem{CakalliSequentialdefinitionsofconnectedness} H. \c{C}akall\i, Sequential definitions of connectedness, Appl. Math. Lett., \textbf{25}(3) (2012), 461-465.
\bibitem{CakalliAVariationonStatisticalWardContinuity} H. Cakalli, A Variation on Statistical Ward Continuity, Bull. Malays. Math. Sci. Soc. DOI 10.1007/s40840-015-0195-0
\bibitem{CakalliAvariationonarithmeticcontinuity} H. Cakalli, {\em A variation on arithmetic continuity}, Bol. Soc. Paran. Mat. \textbf{35}, 3, 195-202, (2017).
\bibitem{CakalliandMucukConnectednessviaasequentialmethod} H. Cakalli and O. Mucuk, On connectedness via a sequential method, Revista de la Union Matematica Argentina, \textbf{54}(2) (2013), 101-109.
\bibitem{CakalliandMucukLacunarystatisticallyupwardanddownwardhalfquasiCauchysequencesJofAnalysis} H. Cakalli and O. Mucuk, {\em Lacunary statistically upward and downward half quasi-Cauchy sequences}, J. Math. Anal. \textbf{7}, 2, 12-23, (2016).
\bibitem{CakalliandTaylanOnabsolutelyalmostconvergence} \c{C}akalli H., and E.Iffet Taylan. On Absolutely Almost Convergence, An. Stiint. Univ. Al. I. Cuza Iasi. Mat. (N.S.),  DOI: 10.2478/aicu-2014-0032 .
\bibitem{ReillyandSubrahmanyamandVamanamurthyCauchysequencesinquasipseudometricspaces} I.L. Reilly, P.V. Subrahmanyam, and M.K. Vamanamurthy, Cauchy sequences in quasipseudometric spaces, Monatsh. Math., \textbf{93}(2) (1982), 127-140.
\bibitem{BurtonandColemanQuasiCauchysequences} D. Burton, J. Coleman, Quasi-Cauchy sequences, Amer. Math. Monthly \textbf{117} (2010) 328-333.
\bibitem{CakalliNthetawardcontinuity} H. Cakalli, N-theta-ward continuity, Abstr. Appl. Anal. \textbf{2012} (2012), Article ID 680456 8 pages.
\bibitem{CakalliSequentialdefinitionsofcompactness} H. \c{C}akalli, Sequential definitions of compactness, Appl. Math. Lett. \textbf{21} (2008), 594-598.
\bibitem{CakalliSlowlyoscillatingcontinuity} H. \c{C}akalli, Slowly oscillating continuity, Abstr. Appl. Anal. \textbf{2008} (2008), Article ID 485706,  5 pages.
\bibitem{CakalliDeltaquasiCauchysequences} H. \c{C}akalli, $\delta$-quasi-Cauchy sequences, Math. Comput. Modelling \textbf{53} (2011), 397-401.
\bibitem{CakalliOnGcontinuity} H. \c{C}akalli, On $G$-continuity, Comput. Math. Appl. \textbf{61} (2011), 313-318.
\bibitem{CakalliStatisticalwardcontinuity} H. \c{C}akalli, Statistical ward continuity. Appl. Math. Lett. \textbf{24} (2011), 1724-1728.
\bibitem{CakalliStatisticalquasiCauchysequences} H. \c{C}akalli, Statistical-quasi-Cauchy sequences, Math. Comput. Modelling 54 (2011) 1620-1624.
\bibitem{CakalliForwardcontinuity} H. \c{C}akalli,  Forward continuity,  J. Comput. Anal. Appl. \textbf{13} (2011), 225-230.
\bibitem{CakalliAvariationonwardcontinuity} H. \c{C}akalli,  A variation on ward continuity, Filomat \textbf{27} 8 (2013), 1545-1549.
\bibitem{CakalliVariationsonquasiCauchysequences} H. \c{C}akall\i, Variations on quasi-Cauchy sequences, Filomat, \textbf{29}(1) (2015), 13-19.
\bibitem{CakalliUpwardanddownwardstatisticalcontinuities} H. \c{C}akall\i,  {\em Upward and downward statistical continuities}, Filomat, \textbf{29}, 10, 2265-2273, (2015).
\bibitem{CakalliOnDeltaquasislowlyoscillatingsequences} H. \c{C}akall\i, {\em On $\Delta$-quasi-slowly oscillating sequences}, Comput. Math. Appl. \textbf{62}, 9, 3567-3574, (2011).
\bibitem{CakalliandAlbayrakNewtypecontinuitiesviaAbelconvergence} H. Cakalli and M. Albayrak, New Type Continuities via Abel Convergence, Scientific World Journal, Volume 2014 (2014), Article ID 398379, 6 pages. http://dx.doi.org/10.1155/2014/398379
\bibitem{CakalliandArasandSonmezLacunarystatisticalwardcontinuity} H. \c{C}akalli, C.G. Aras, A. Sonmez, Lacunary statistical ward continuity, AIP Conf. Proc. \textbf{1676}, 020042 (2015); http://dx.doi.org/10.1063/1.4930468
\bibitem{CakalliandCanakPnsabsolutealmostconvergentsequences} H.\c{C}akalli, and G.\c{C}anak, $(P_{n} ,s)$-absolute almost convergent sequences, Indian J. Pure Appl. Math., 28 4 (1997) 525-532.
\bibitem{CanakandDik} I. Canak and M. Dik, New types of continuities, Abstr. Appl. Anal. \textbf{201}0 (2010), Article ID 258980, 6 pages.
\bibitem{CakalliandHazarikaIdealquasiCauchysequences} H. \c{C}akalli,  and B. Hazarika, Ideal quasi-Cauchy sequences, J. Inequal. Appl. \textbf{2012} (2012), Article 234,  11 pages.
\bibitem{CakalliandKaplanAstudyonNthetaquasiCauchysequences} H. \c{C}akalli, and H. Kaplan, A study on N-theta-quasi-Cauchy sequences,  Abstr. Appl. Anal.  2013 (2013), Article ID 836970, 4 pages.
\bibitem{CakalliandKaplanAvariationonlacunarystatisticalquasiCauchysequencesCommunications} H. Cakalli, and H. Kaplan, {\em A variation on lacunary statistical quasi Cauchy sequences}, Communications Faculty of Sciences University of Ankara-Series A1 Mathematics and Statistics, \textbf{66}, 2, 71-79, (2017). 10.1501/Commua1 0000000802
\bibitem{CakalliandPattersonFunctionspreservingslowlyoscillatingdoublesequences} H. \c{C}akalli and R.F. Patterson, Functions preserving slowly oscillating double sequences, An. Stiint. Univ. Al. I. Cuza Iasi. Mat. (N.S.) in press.
\bibitem{CakalliandDasFuzzycompactnessviasummability} H. \c{C}akalli, and Pratulananda Das, Fuzzy compactness via summability, Appl. Math. Lett. \textbf{22} (2009), 1665-1669.
\bibitem{CakalliandSonmezSlowlyoscillatingcontinuityinabstractmetricspaces} H. \c{C}akalli and A. Sonmez, Slowly oscillating continuity in abstract metric spaces, Filomat, \textbf{27} (2013), 925-930.
\bibitem{CakalliandSonmezandArasLambdastatisticallywardcontinuity} H. \c{C}akalli, A. Sonmez, and C.G. Aras, $\lambda$-statistical ward continuity, An. Stiint. Univ. Al. I. Cuza Iasi. Mat. (N.S.) DOI: 10.1515/aicu-2015-0016 March 2015.
\bibitem{CakalliandSonmezandGenc}  H. \c{C}akall\i, A. Sonmez, and C. Genc, On an equivalence of topological vector space valued cone metric spaces and metric spaces, Appl. Math. Lett. \textbf{25} (2012), 429-433.
\bibitem{DjurcicandKocinacandZizovicDoublesequencesandselections} Djurcic, D.,  Kocinac, Ljubisa D. R., Zizovic, M. R. Double Sequences and Selections, Abstr. Appl. Anal. Hindawi Publ. Corp., New York, Article Number: 497594, 2012.  DOI: 10.1155/2012/497594
\bibitem{ErdemandArasandCakalliandSonmezSoftmatricesonsoftmultisetsinanoptimaldecisionprocess} A.E. Coskun, C.G Aras, H. Cakalli, and A. Sonmez, {\em Soft matrices on soft multisets in an optimal decision process}, AIP Conference Proceedings, \textbf{1759}, 1,  020099 (2016); doi: 10.1063/1.4959713
\bibitem{Fast} H. Fast, Sur la convergence statistique, Colloq. Math.  \textbf{2} (1951), 241-244.
\bibitem{FridyOnstatisticalconvergence} J.A. Fridy, On statistical convergence, Analysis \textbf{5} (1985), 301-313.
\bibitem{KaplanandCakalliVariationsonstronglylacunaryquasiCauchysequencesAIP} H. Kaplan, H. Cakalli, {\em Variations on strongly lacunary quasi Cauchy sequences}, AIP Conf. Proc. \textbf{1759}, Article Number: 020051, (2016). doi: http://dx.doi.org/10.1063/1.4959665
\bibitem{KaplanandCakalliVariationsonstronglacunaryquasiCauchysequencesJNonlinearSciAppl} H. Kaplan, H. Cakalli, {\em Variations on strong lacunary quasi-Cauchy sequences}, J. Nonlinear Sci. Appl. \textbf{9}, 4371-4380, (2016).
\bibitem{KaplanandCakalliVariationsonstronglylacunaryquasiCauchysequencesAIP} H. Kaplan, H. Cakalli, Variations on strongly lacunary quasi Cauchy sequences, AIP Conf. Proc. \textbf{1759} (2016) Article Number: 020051


\bibitem{KocinacSelectionpropertiesinfuzzymetricspaces} Ljubisa D.R. Ko\v cinac, Selection properties in fuzzy metric spaces, Filomat,  \textbf{26}(2) (2012), 99-106.
\bibitem{MaioKocinac} G. Di Maio, and Lj.D.R. Ko\v cinac, Statistical convergence in topology, Topology Appl. \textbf{156} (2008), 28-45.
\bibitem{MucukandCakalliGsequentiallyconnectednessfortopologicalgroupswithoperations} O. Mucuk and H. Cakalli, {\em G-sequentially connectedness  for topological groups with operations}, AIP Conference Proceedings, \textbf{1759}, 020038, (2016). doi: 10.1063/1.4959652
\bibitem{MucukandSahanOnGsequentialcontinuity} O. Mucuk, T. \c{S}ahan¸ On $G$-Sequential Continuity, Filomat, \textbf{28}(6) (2014), 1181-1189.
\bibitem{MursaleenLamdastatisticalconvergence} M. Mursaleen, $\lambda$-statistical convergence, Math. Slovaca \textbf{50} (2000), 111-115.
\bibitem{PalandSavasandCakalliIconvergenceonconemetricspaces}	S.K. Pal,  E. Savas, and H. Cakalli, $I$-convergence on cone metric spaces, Sarajevo J. Math. \textbf{9}   (2013), 85-93.
\bibitem{PalladinoOnhalfcauchysequences} F.J. Palladino,  On half Cauchy sequences, Arxiv. \textit{arXiv}:1102.4641v1, (2012), 3 pages.
\bibitem{PattersonandCakalliQuasiCauchydoublesequences} R.F. Patterson and H. Cakalli, Quasi Cauchy double sequences, Tbilisi Mathematical Journal \textbf{8}(2)  (2015), 211-219.
\bibitem{PattersonandSavasRateofPconvergenceoverequivalenceclassesofdoublesequencespaces} R.F. Patterson, and E. Sava\c{s}, Rate of P-convergence over equivalence classes of double sequence spaces, Positivity  \textbf{16}(4) (2012), 739-749.
\bibitem{PattersonandSavasAsymptoticEquivalenceofDoubleSequences} R.F. Patterson, and E. Savas, Asymptotic equivalence of double sequences, Hacet. J. Math. Stat.  \textbf{41} (2012), 487-497.

\bibitem{SavasandCakalliIdealstatisticalquasiCauchysequences} E. Savas, H. \c{C}akall\i, {\em Ideal statistically quasi Cauchy sequences}, AIP Conf. Proc. \textbf{1759}, 020057, (2016). doi: 10.1063/1.4959671
\bibitem{SonmezParacompactness} Sonmez, A., On paracompactness in cone metric spaces, Appl. Math. Lett. \textbf{23}, (2010), 494-497.
\bibitem{SonmezandCakalliConenormedspacesandweightedmeans} A. Sonmez, and H. \c{C}akall\i, Cone normed spaces and weighted means, Math. Comput. Modelling  \textbf{52}  (2010), 1660-16660.
\bibitem{VallinCreatingslowlyoscillatingsequencesandslowlyoscillatingcontinuousfunctions} R.W. Vallin, Creating slowly oscillating sequences and slowly oscillating continuous functions (with an appendix by Vallin and H. \c{C}akalli), Acta Math. Univ. Comenianae \textbf{25} (2011), 71-78.
\bibitem{YayingandHazarikandCakalliNewresultsinquasiconemetricspaces} T. Yaying, B. Hazarika, H. Cakalli, {\em New results in quasi cone metric spaces}, J. Math. Computer Sci. \textbf{16}, 435-444, (2016).
\bibitem{YildizIstatistikselboslukludelta2quasiCauchydizileri} \c{S}. Y\i ld\i z, {\em \.{I}statistiksel bo\c{s}luklu delta 2 quasi Cauchy dizileri}, Sakarya University Journal of Science, \textbf{21}, 6, (2017).
\end{thebibliography}
\end{document}